\documentclass[submission]{FPSAC2024}

\usepackage{pdflscape}  
\usepackage[utf8]{inputenc}
\usepackage[T1]{fontenc}
\usepackage{amsfonts, wasysym}
\usepackage{graphics, graphicx}
\usepackage{multido, xargs, colortbl, multicol, multirow}
\usepackage{calc, ifthen, xstring, stmaryrd, blkarray, stackrel}
\usepackage{listings, hvfloat}
\usepackage[normalem]{ulem}
\usepackage{url, hyperref, hypcap}
\hypersetup{colorlinks=true, citecolor=blue, linkcolor=blue}

\usepackage[skip=0.333\baselineskip]{caption}
\usepackage{subcaption}

\usepackage{pgfkeys}
\usepackage{tikz}
\usetikzlibrary{trees, shapes, arrows, arrows.meta, matrix, calc, graphs,
  plotmarks, backgrounds}
\usetikzlibrary{decorations,
  decorations.markings,decorations.pathmorphing,decorations.pathreplacing}
\pgfdeclaredecoration{cappedcurveto}{initial}{%
  \state{initial}[width=\pgfdecoratedinputsegmentlength/100]
  {\pgfpathlineto{\pgfpointorigin} }%
  \state{final}{}%
}%
\tikzset{
  draw left/.style={
    decorate, decoration={cappedcurveto, raise={.5*\pgflinewidth}}
  },
  draw right/.style={
    decorate, decoration={cappedcurveto, raise={-.5*\pgflinewidth}}
  }
}

\tikzset{blueNode/.style={blue, edge from parent/.style={blue, draw}}}
\tikzset{redNode/.style={red, edge from parent/.style={red, draw}}}
\tikzset{greenNode/.style={green, edge from parent/.style={green, draw}}}
\tikzset{cyanNode/.style={cyan, edge from parent/.style={cyan, draw}}}
\tikzset{blueEdge/.style={edge from parent/.style={blue, draw}}}
\tikzset{redEdge/.style={edge from parent/.style={red, draw}}}
\tikzset{violetEdge/.style={edge from parent/.style={violet, draw}}}

\usetikzlibrary{cd}
\AtBeginEnvironment{tikzcd}{\tikzexternaldisable}
\AtEndEnvironment{tikzcd}{\tikzexternalenable}

\makeatletter
\newcommand{\myitem}[1]{%
\item[#1]\protected@edef\@currentlabel{#1}%
}
\makeatother

\crefname{enumi}{Condition}{Conditions}

\newcommand{\lststyle}{\upshape\ttfamily}

\title[Diagrammatic Okada Algebra]{%
  Diagram model for the Okada algebra and monoid }

\author[F.~Hivert \and J.~Scott]{Florent Hivert\thanks{\href{mailto:Florent.Hivert@universite-paris-saclay.fr}{Florent.Hivert@universite-paris-saclay.fr}}\addressmark{1}, \and Jeanne Scott\thanks{\href{mailto:jeanne@imsc.res.in}{jeanne@imsc.res.in}}\addressmark{2}}

\address{\addressmark{1}Laboratoire Interdisciplinaire des Sciences du
  Num\'erique (UMR CNRS 9015) \\
  1, rue Raimond Castaing, 91190 GIF-SUR-YVETTE, FRANCE \\
  Université Paris-Saclay, CNRS \\ 
  \addressmark{2}Department of Mathematics, Brandeis University \\
  415 South Street, Waltham, MA 02453 \\ } 

\received{\today}

\newtheorem{theorem}{Theorem}[section]

\newtheorem{proposition}[theorem]{Proposition}
\newtheorem{lemma}[theorem]{Lemma}
\newtheorem{definition}[theorem]{Definition}

\theoremstyle{definition}

\newtheorem{remark}[theorem]{Remark}

\newcommand{\K}{\mathbb{K}} 
\newcommand{\C}{\mathbb{C}} 

\newcommand{\set}[2]{\left\{ #1 \;\middle|\; #2 \right\}} 
\newcommand{\eqdef}{\mbox{\,\raisebox{0.2ex}{\scriptsize\ensuremath{\mathrm:}}\ensuremath{=}\,}}
\newcommand{\card}{\operatorname{\#}}            

\newcommand{\JJ}{\mathrel{\mathcal{J}}}

\newcommand{\RR}{\mathrel{\mathcal{R}}}

\newcommand{\ie}{\textit{i.e.}~} 
\newcommand{\qandq}{\quad\text{and}\quad} 
\newcommand{\defn}[1]{\emph{\textsf{\color{blue} #1}}} 

\definecolor{green}{RGB}{57,181,74} 

\newcommand{\bl}{|[text=white,fill=black]|}%
\newcommand{\nn}{|[draw=none]|}%
\tikzset{diamond/node/.style={
    scale=0.70, minimum size=0.9cm, inner sep=0mm,
    draw, regular polygon, regular polygon sides=4,
    yscale=-1, rotate=45, shape border rotate=45
  }}%
\tikzset{diamond/matrix/.style={matrix of nodes,
  transform canvas={yscale=-1, rotate=45},
  nodes={diamond/node},row sep=-\pgflinewidth,column sep=-\pgflinewidth}}
\tikzset{diamond/.style={
    baseline={([yshift=-.5ex]current bounding box.center)},
  }}
\pgfkeys{/tikz/diamond/boundabs/.code 2 args={
  \useasboundingbox[scale=0.5] ($(#1)$) rectangle ($(#2)$);
}}
\pgfkeys{/tikz/diamond/bound/.code 2 args={
    \useasboundingbox[scale=0.5]
    ($(-0.1, -0.9) + #1*(-0.62, 0) + #2*(0, 1.24)$)
    rectangle
    ($(0.2, -0.4) + #1*(0.62, 0)$);
}}
\tikzset{debug/.style={
    show background rectangle,    
  }}
\newcommand{\stC}[2]{
  \node[diamond]{};
  \begin{scope}
    \clip (-6.2pt,-6.2pt) rectangle (6.2pt,6.2pt);
    \draw[line width=2pt, #1] (-7.32pt,1pt) -- (1pt,-7.32pt);
    \draw[line width=2pt, #2] (7.32pt,-1pt) -- (-1pt,7.32pt);
  \end{scope}
}
\newcommand{\st}{\stC{black}{black}}
\newcommand{\tuC}[2]{
  \node[diamond]{};
  \begin{scope}
    \clip (-6.2pt,-6.2pt) rectangle (6.2pt,6.2pt);
    \draw[line width=2.4pt, #1] (-4.5pt,4.5pt) circle (4.8pt);
    \draw[line width=2.4pt, #2] (4.5pt,-4.5pt) circle (4.8pt);
  \end{scope}
}
\newcommand{\tu}{\tuC{black}{black}}
\newcommand{\stbC}[1]{
  \node[diamond, draw=none]{};
  \begin{scope}
    \clip (-6.2pt,-6.2pt) rectangle (6.2pt,6.2pt);
    \draw[line width=2pt, #1] (7.32pt,-1pt) -- (-1pt,7.32pt);
  \end{scope}
}
\newcommand{\stb}{\stbC{black}}
\newcommand{\sttC}[1]{
  \node[diamond, draw=none]{};
  \begin{scope}
    \clip (-6.2pt,-6.2pt) rectangle (6.2pt,6.2pt);
    \draw[line width=2pt, #1] (-7.32pt,1pt) -- (1pt,-7.32pt);
  \end{scope}
}
\newcommand{\stt}{\sttC{black}}
\newcommand{\TLlab}[1]{\nn\raisebox{-8pt}[0pt][0pt]{$#1$}}%

\newcommand{\ov}[1]{\overline{#1}}

\newcommand{\arctip}{\rule[0.1ex]{0.5pt}{1ex}}
\newcommand{\arcdrawing}[1]{\arctip\rule[0.55ex]{#1}{0.5pt}}%

\newlength{\arcwidth}
\newlength{\arcwidthB}
\newsavebox{\arclabel}
\newsavebox{\arcdiag}
\newcommand{\arc}[4][]{%
  \savebox{\arclabel}{$\scriptstyle#3$}%
  \savebox{\arcdiag}{$\scriptstyle#1$}%
  \settowidth{\arcwidth}{\usebox{\arclabel}}%
  \settowidth{\arcwidthB}{\usebox{\arcdiag}}%
  \ifthenelse{\lengthtest{\arcwidthB>\arcwidth}}{%
    \setlength{\arcwidth}{\arcwidthB}}{}%
  \addtolength{\arcwidth}{1em}%
  \ifthenelse{\equal{#4}{}}{\def\endtip{}}{\def\endtip{\arctip}}%
  \ifthenelse{\equal{#1}{}}{%
    #2\mathrel{\stackrel{\usebox{\arclabel}}{\arcdrawing{\arcwidth}\endtip}}#4%
  }{#2\mathrel{\stackrel[{\usebox{\arcdiag}}]{\usebox{\arclabel}}%
      {\arcdrawing{\arcwidth}\endtip}}#4%
  }%
}
\newcommand{\harc}[3][{}]{\arc[#1]{#2}{#3}{}}

\newcommand{\glue}[2]{#1\Join#2}  
\newcommand{\leftd}[1]{\langle #1|} 
\newcommand{\rightd}[1]{|#1\rangle} 

\newcommand{\isdomeq}{\preceq}  
\newcommand{\isdom}{\prec}  
\newcommand{\emptyword}{\varepsilon}

\newcommandx{\setN}[1][1=N]{[#1]}
\newcommandx{\setNN}[1][1=N]{\overline{[#1]}}
\newcommandx{\SG}[1][1=N]{\mathfrak{S}_{#1}} 
\newcommandx{\Okada}[1][1=N]{\mathsf{O}_{#1}}
\newcommandx{\OkArc}[1][1=N]{\widetilde{\mathsf{O}}_{#1}}
\newcommandx{\TL}[1][1=N]{\mathsf{TL}_{#1}}
\newcommand{\E}[1]{\mathsf{E}_{#1}}

\newcommand{\G}[1]{\mathsf{G}_{#1}}

\newcommand{\id}{\mathsf{id}}      
\newcommand{\code}{\mathsf{code}}  
\newcommandx{\lab}[1][1=]{\operatorname{lab}_{#1}}

\newcommand{\OSys}{\mathsf{O}}

\newcommand{\YFS}{\mathbb{YFS}}
\newcommandx{\YFSn}[1][1=N]{\YFS(#1)}
\newcommand{\YF}{\mathbb{YF}}
\newcommandx{\YFn}[1][1=N]{\YF(#1)}

\DeclareMathOperator{\PropLab}{PLab}   
\DeclareMathOperator{\Chain}{Chain}

\DeclareMathOperator{\RS}{RS}

\abstract{
  
  It is well known that the Young lattice is the Bratelli diagram of the
  symmetric groups expressing how irreducible representations restrict from
  $\SG[N]$ to $\SG[N-1]$. In 1988, Stanley discovered a similar lattice called
  the Young-Fibonacci lattice which was realized as the Bratelli diagram of a
  family of algebras by Okada in 1994.

  In this paper, we realize the Okada algebra and its associated monoid
  using a labeled version of Temperley-Lieb arc-diagrams. 
  We prove in full generality that the dimension of
  the Okada algebra is $n!$. In particular, we interpret a natural bijection
  between permutations and labeled arc-diagrams as an instance of 
  Fomin's Robinson-Schensted correspondence for the Young-Fibonacci
  lattice. We prove that the Okada monoid is aperiodic and describe its Green
  relations. Lifting those results to the algebra allows us to construct a cellular
  basis of the Okada algebra.  }

\resume{
  Il est bien connu que le treillis de Young peut s'interpréter comme le
  diagramme de Bratelli des groupes symétriques, décrivant, par exemple,
  comment les représentations irréductibles se restreignent de $\SG[n]$ à
  $\SG[n-1]$. En 1975, Stanley a découvert un treillis similaire appelée
  treillis de Young-Fibonacci qui a été interprété comme le diagramme de
  Bratelli d'une famille d'algèbres par Okada en 1994.

  Dans cet article, nous réalisons l'algèbre d'Okada et le monoïde associé
  grâce à une version étiquetée des diagrammes d'arcs du monoïde de Jones et
  de l'algèbre de Tempeley-Lieb. Nous prouvons en toute généralité que
  l'algèbre d'Okada est de dimension $n!$. En particulier, nous interprétons
  la bijection naturelle entre les permutations et les diagrammes d'arcs comme
  une instance de la correspondance de Robinson-Schensted-Fomin associée au
  treillis de Young-Fibonacci. Nous prouvons que le monoïde est apériodique et
  décrivons ses relations de Green. En relevant, ces dernières à l'algèbre
  nous en construisant une base cellulaire.  }

\begin{document}
\maketitle

\section{Introduction}

The theory of $1$-differential posets was developed by
R. Stanley~\cite{Stanley1988} as a framework for generalizing the
Robinson-Schensted correspondence beyond the combinatorics of the Young
lattice $\mathbb{Y}$ of integer partitions.  A similar undertaking was made by
S. Fomin in his work on dual graded graphs and growth processes, where the later
technique was used to construct an explicit RS-correspondence for
Stanley's \defn{Young-Fibonacci} lattice $\YF$~\cite{Fomin1995,Stanley1988}.
Both $\mathbb{Y}$ and $\YF$ are $1$-differential and they
are the only lattices having this property.
Fomin's approach involves a Fibonacci version of standard tableaux; a notion
later examined independently by T. Roby, K. Killpatrick, and J. Nzeutchap
(whose formulation by-passes the growth construction altogether),
see~\cite{Nzeutchap2009} and the references therein.

S. Okada~\cite{Okada1994} showed that the $\YF$-lattice supports a theory
of {\it clone symmetric functions} with analogues of the classical bases
(e.g. complete homogeneous, Schur, and power-sum symmetric functions) as well as a
$\YF$-variant of the Littlewood-Richardson rule.  The clone theory
appears in Goodman-Kerov's determination of the {\it Martin boundary} of the
$\YF$-lattice~\cite{Goodman1997} and is related to various random
processes.

The Okada algebras $\{ \Okada(X,Y) \}_{\scriptscriptstyle \mathrm{N} \geq 0}$
were introduced by S. Okada as a counterpart to the clone theory, 
and occupy a role similar to that
played by the symmetric groups in the classical theory
of symmetric functions. Okada
algebras are finite dimensional, associative, and depend on
parameters $X=(x_1,\dots,x_{{\scriptscriptstyle N}-1})$ and
$Y=(y_1,\dots, y_{{\scriptscriptstyle N}-2})$. When those parameters are
generic, they are semi-simple and their branching rule, which
describes how irreducible representations restrict from
$\Okada[\scriptscriptstyle N](X,Y)$ to
$\Okada[{\scriptscriptstyle N}-1](X,Y)$, is expressed by the covering
relations of the $\YF$-lattice.



In this paper we  
realize the Okada algebra $\Okada(X,Y)$
as a diagram algebra
with a multiplicative/monoidal basis
expressed in terms of certain arc-labeled, non-crossing
perfect matchings (as appear in both  the Temperley-Lieb
and Martin-Saleur Blob algebras~\cite{Martin1994a}).
Like most diagram algebras, this basis is cellular and
affords us with a 
novel, diagrammatic presentation of the irreducible 
representations of $\Okada(X,Y)$
(i.e. as \defn{cell modules}). 
We interpret Fomin's RS-correspondence
diagrammatically. This involves constructing a 
bijection between saturated chains in the $\YF$-lattice
(presented in terms of sequences of \defn{Fibonacci sets})
and \defn{Okada half arc-diagrams}. 
In addition we examine the structure theory of the Okada algebra
and monoid via a \defn{dominance order} on Fibonacci sets.

\paragraph{Aknowledgment:} F. Hivert would like to thank James
Mitchell and Nicolas Thiéry for fruitful discussions. Likewise, J. Scott 
thanks Philippe Biane and  
Anatoly Vershik for their input and guidance.
The computations were
done using the SageMath computer algebra system.

\section{Background}
Throughout this paper $N$ denotes a non-negative integer. We denote by $\setN$
the set $\{1,\dots,N\}$. We often write negative numbers as overlined numbers
such as $\ov4$. The cardinality of a set $S$ is denoted $\card S$. For a
non-negative integer $N$ we endow $[N] \cup [\ov{N}]$ with the total order
$
\{1 < 2< \dots < N < \ov{N} < \dots < \ov{2} < \ov{1}\}\,.
$
Overlining numbers which are negative should also
help the reader remember this unusual ordering.

Stanley's original construction of the Young-Fibonacci lattice~\cite{Stanley1988}
involves endowing the set of \defn{Fibonacci words}, i.e. binary words in the
alphabet $\{1,2\}$, with a partial order.
We present an alternative
description using \defn{Fibonacci sets}.

\begin{definition}
 A \defn{Fibonacci set of rank $N$} is a subset $S = \{s_1<s_2<\dots<s_k\}$ of
  $\setN$ whose size $k$ has the same parity as $N$ and such that $s_\ell$ have
  the same parity as $\ell$. We write $\YFS_N$ for the
  collection of all rank $N$ Fibonacci sets and $\YFS$ for the
  \emph{disjoint} union of $\YFS_N$ as $N$ varies.
\end{definition}
The entire interval $\setN$ itself is always a Fibonacci set of rank $N$,
while $\emptyset$ is a Fibonacci set only when $N$ is even. We emphasize on
the fact that in $\YFS$ the set $\{1,2,5\}$ of rank $5$ is not the same
Fibonacci set as $\{1,2,5\}$ of rank $7$. When they need to be distinguished
we include $N$ as a subscript, as in $\{1,2,5\}_5$ and $\{1,2,5\}_7$.

The covering relations which generate the lattice structure on $\YFS$ are
defined by $S \lhd T$ if and only if
$S \in \YFS_{\scriptscriptstyle \mathrm{N-1}}$ and
$T \in \YFS_{\scriptscriptstyle \mathrm{N}}$ and one of these two sets
can be obtained from the other one by removing its largest element. 
Stanley's description is equivalent to ours through the bijection 
sending a binary word $w$ to the set of the sums of its suffixes
whose first digit is $1$.
The Hasse diagram of $\YFS$ upto rank~$5$ is illustrated below.
\bigskip

 
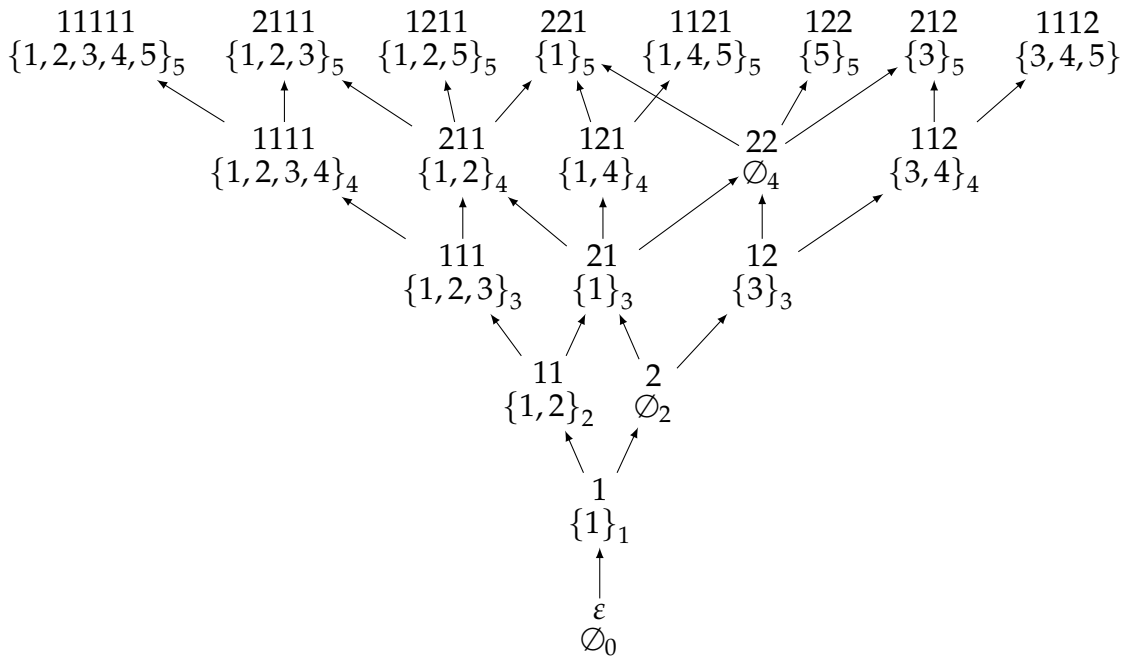
\begin{tikzpicture}[>=latex,xscale=1.2,yscale=0.9,line join=bevel]
\tikzset{every node/.style={rectangle split,rectangle split parts=2,inner sep=0.5mm}}
\node (node_0) at (163.0bp,7.5bp) {$\emptyword$ \nodepart{second} $\emptyset_0$};
 \node (node_1) at (163.0bp,57.5bp) {$1$ \nodepart{second} $\left\{1\right\}_1$};
 \node (node_2) at (147.0bp,106.5bp) {$11$ \nodepart{second} $\left\{1, 2\right\}_2$};
 \node (node_6) at (180.0bp,106.5bp) {$2$ \nodepart{second} $\emptyset_2$};
 \node (node_3) at (120.0bp,155.5bp) {$111$ \nodepart{second} $\left\{1, 2, 3\right\}_3$};
 \node (node_10) at (164.0bp,155.5bp) {$21$ \nodepart{second} $\left\{1\right\}_3$};
 \node (node_4) at (64.0bp,204.5bp) {$1111$ \nodepart{second} $\left\{1, 2, 3, 4\right\}_4$};
 \node (node_13) at (120.0bp,204.5bp) {$211$ \nodepart{second} $\left\{1, 2\right\}_4$};
 \node (node_5) at (5.0bp,253.5bp) {$11111$ \nodepart{second} $\left\{1, 2, 3, 4, 5\right\}_5$};
 \node (node_15) at (64.0bp,253.5bp) {$2111$ \nodepart{second} $\left\{1, 2, 3\right\}_5$};
 \node (node_7) at (214.0bp,155.5bp) {$12$ \nodepart{second} $\left\{3\right\}_3$};
 \node (node_8) at (268.0bp,204.5bp) {$112$ \nodepart{second} $\left\{3, 4\right\}_4$};
 \node (node_16) at (214.0bp,204.5bp) {$22$ \nodepart{second} $\emptyset_4$};
 \node (node_9) at (310.0bp,253.5bp) {$1112$ \nodepart{second} $\left\{3, 4, 5\right\}$};
 \node (node_18) at (268.0bp,253.5bp) {$212$ \nodepart{second} $\left\{3\right\}_5$};
 \node (node_11) at (164.0bp,204.5bp) {$121$ \nodepart{second} $\left\{1, 4\right\}_4$};
 \node (node_12) at (195.0bp,253.5bp) {$1121$ \nodepart{second} $\left\{1, 4, 5\right\}_5$};
 \node (node_19) at (152.0bp,253.5bp) {$221$ \nodepart{second} $\left\{1\right\}_5$};
 \node (node_14) at (112.0bp,253.5bp) {$1211$ \nodepart{second} $\left\{1, 2, 5\right\}_5$};
 \node (node_17) at (235.0bp,253.5bp) {$122$ \nodepart{second} $\left\{5\right\}_5$};
 \draw [black,->] (node_0) -- (node_1);
 \draw [black,->] (node_1) -- (node_2);
 \draw [black,->] (node_1) -- (node_6);
 \draw [black,->] (node_2) -- (node_3);
 \draw [black,->] (node_2) -- (node_10);
 \draw [black,->] (node_3) -- (node_4);
 \draw [black,->] (node_3) -- (node_13);
 \draw [black,->] (node_4) -- (node_5);
 \draw [black,->] (node_4) -- (node_15);
 \draw [black,->] (node_6) -- (node_7);
 \draw [black,->] (node_6) -- (node_10);
 \draw [black,->] (node_7) -- (node_8);
 \draw [black,->] (node_7) -- (node_16);
 \draw [black,->] (node_8) -- (node_9);
 \draw [black,->] (node_8) -- (node_18);
 \draw [black,->] (node_10) -- (node_11);
 \draw [black,->] (node_10) -- (node_13);
 \draw [black,->] (node_10) -- (node_16);
 \draw [black,->] (node_11) -- (node_12);
 \draw [black,->] (node_11) -- (node_19);
 \draw [black,->] (node_13) -- (node_14);
 \draw [black,->] (node_13) -- (node_15);
 \draw [black,->] (node_13) -- (node_19);
 \draw [black,->] (node_16) -- (node_17);
 \draw [black,->] (node_16) -- (node_18);
 \draw [black,->] (node_16) -- (node_19);
\end{tikzpicture}

\begin{definition}
Fix a positive integer $N$. Given a field $\K$, fix also $X=(x_1,\dots,x_{N-1})$
and $Y=(y_1,\dots, y_{N-2})$ two sequences of parameters in $\K$.
The \defn{Okaka algebra $\Okada(X, Y)$} is the algebra generated by
$\set{\E{i}}{i=1\dots{N-1}}$ with the relations
\begin{alignat}{3}
  \E{i}^2&=x_i\E{i} &\qquad& 1\leq i \leq N-1, \tag{$I(X,Y)$}\label{relIxy}\\
  \E{i}\E{j}&=\E{j}\E{i} &\qquad&\vert i-j\vert \geq 2,   \tag{C(X,Y)}\label{relCxy}\\
  \E{i+1}\E{i}\E{i+1}&=y_i\E{i+1} &\qquad&  1\leq i \leq N-2, \tag{S(X,Y)}\label{relSxy}
\end{alignat}
\end{definition}
If all the $X$'s and the $Y$'s are equal to $1$, the Okada algebra is actually
the algebra of a monoid; we call this the \defn{Okaka Monoid} and
denote it $\Okada$. Recall that setting all $y_i\eqdef 1$ and all
$x_i\eqdef q$ and adding the extra relation $\E{i}\E{i+1}\E{i+1}=\E{i}$
defines the Temperley-Lieb algebra which is also a deformation of the algebra
of a monoid called the Jones monoid (obtained when $q=1$). \medskip

We now review some of Okada's results~\cite{Okada1994}:
For generic values of the $X$ and $Y$ parameters
$\Okada(X, Y)$ is semi-simple and its
irreducible representations $V_T$
correspond to rank $N$ Fibonacci sets $T$. When $V_T$ is restricted to
the subalgebra $\Okada[N-1](X,Y) \subset \Okada(X,Y)$ 
it decomposes as a direct sum of irreducible representations
$V_S$ of $\Okada[N-1](X,Y)$ where $S \lhd T$
is a covering relation in $\YFS$.

The dimension of $\Okada(X,Y)$ is $N!$ and a basis for $\Okada(X,Y)$ can be
constructed from permutations in the following way.  Recall that the
\defn{code} of a permutation $\sigma\in\SG[N]$ is
$\code(\sigma)=(c_1,\dots c_N)$ where
$c_i \eqdef \card\set{j<i}{\sigma^{-1}(j)>\sigma^{-1}(i)}$. It is well known that the
product $\prod_{i=1}^n \sigma_{i-1}\sigma_{i-2} \cdots \sigma_{i - c_i}$ taken from left to
right, increasing with $i$, is the lexicographically minimal reduced factorization of
$\sigma$ into simple transpositions $\sigma_i =(i, i+1)$. Define
$\E{\sigma}\eqdef\prod_{i=1}^n\E{i-1}\E{i-2} \cdots \E{i-c_i}$. Okada showed
in~\cite{Okada1994} that the family $\{\E{\sigma} \! \mid \! \sigma\in\SG[n]\}$ is,
generically, a basis of the Okada algebra. His proof, however, requires
semi-simplicity and doesn't apply to degenerate specializations, such as the
monoid case.


\section{Diagram models for the Okada Monoid and Algebra}

The goal of this section is to build a basis of the Okada algebra in full
generality using rewriting techniques. Inspired by Viennot's theory of heaps of
dimers~\cite{ViennotHeap}, we use diagram rewriting rather than
word rewriting.

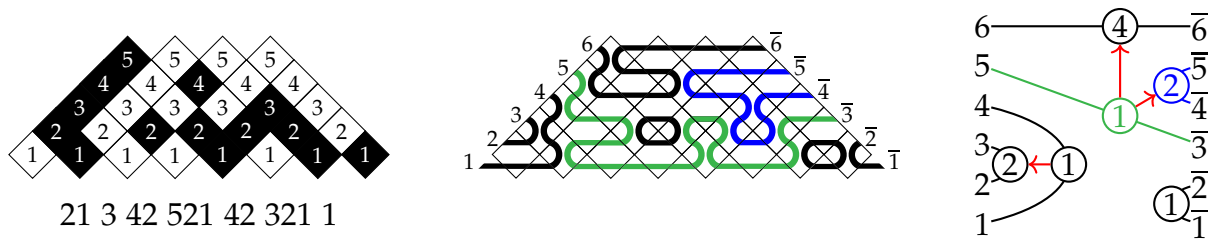
\begin{figure}[ht]
  \vskip-0.5cm
  \[
\begin{tikzpicture}[diamond, diamond/bound={8}{3}]\matrix[diamond/matrix] {
    1  \\
  \bl2&\bl1\\
  \bl3&  2 &  1 \\
  \bl4&  3 &\bl2&  1 \\
  \bl5&  4 &  3 &\bl2&\bl1\\
  &  5 &\bl4&  3 &\bl2&  1 \\
  &   &  5 &  4 &\bl3&\bl2&\bl1\\
  &   &   &  5 &  4 &  3 &  2 &\bl1\\
};
\node at (0,-0.8) {$21\ 3\ 42\ 521\ 42\ 321\ 1$};
\end{tikzpicture}
\hspace{1cm}
\begin{tikzpicture}[diamond, diamond/bound={9}{3}]\matrix[diamond/matrix] {
  \TLlab1 & \stt\\
  \TLlab2 &\st& \stt\\
  \TLlab3 &\tu&\tuC{black}{green}& \sttC{green}\\
  \TLlab4 &\tuC{black}{green}&\stC{green}{green}&\stC{green}{green}& \sttC{green}\\
  \TLlab5 &\tuC{green}{black}&\stC{black}{green}&\tuC{green}{black}&\stC{black}{green}& \sttC{green}\\
  \TLlab6 &\tu&\st&\st&\tuC{black}{green}&\tuC{green}{green}& \sttC{green}\\
          &\stb&\st&\tuC{black}{blue}&\st\stC{blue}{green}&\tuC{green}{blue}&\stC{blue}{green}& \sttC{green}\\
          &   &\stb&\stC{black}{blue}&\stC{blue}{blue}&\tuC{blue}{blue}&\tuC{blue}{green}&\tuC{green}{black}& \stt\\
          &   &   &\stb&\stC{black}{blue}&\stC{blue}{blue}&\stC{blue}{green}&\stC{green}{black}&\tu& \stt\\
  & & & & \TLlab{\ov6} & \TLlab{\ov5} &\TLlab{\ov4} &\TLlab{\ov3} &\TLlab{\ov2} &\TLlab{\ov1} \\
};
\end{tikzpicture}
\hspace{1cm}
\begin{tikzpicture}[yscale=0.35,xscale=0.65, thick, baseline={(current bounding box.east)}]
\tikzstyle{vertex} = [shape=rectangle, minimum height=15pt, minimum width=0pt, inner sep=0pt]
\tikzstyle{mid} = [draw, fill=white, shape=circle, minimum size=2pt, inner sep=1pt]
\foreach \i in {1,...,6} {
  \node[vertex] (G-\i) at ($(-3, \i*1.5-1.5)$) {$\i$};
  \node[vertex] (G--\i) at ($(1.44, \i*1.5-1.5)$) {$\ov{\i}$};
}
\draw[green] (G--3) -- (G-5) node[pos=0.35294117647058826,mid] (M-3-5) { 1 };
\draw (G-1) .. controls +(2.3, 1.05) and +(2.3, -1.05) .. (G-4) node[pos=0.5,mid] (M1-4) { 1 };
\draw[blue] (G--5) .. controls +(-0.7, -0.35) and +(-0.7, 0.35) .. (G--4) node[pos=0.5,mid] (M-5--4) { 2 };
\draw (G--6) -- (G-6) node[pos=0.35294117647058826,mid] (M-6-6) { 4 };
\draw (G-2) .. controls +(0.7, 0.35) and +(0.7, -0.35) .. (G-3) node[pos=0.5,mid] (M2-3) { 2 };
\draw (G--2) .. controls +(-0.7, -0.35) and +(-0.7, 0.35) .. (G--1) node[pos=0.5,mid] (M-2--1) { 1 };
\draw[color=red, <-] (M2-3) -- (M1-4);
\draw[color=red, <-] (M-5--4) -- (M-3-5);
\draw[color=red, <-] (M-6-6) -- (M-3-5);
\node at (0,8.5) {};    
\end{tikzpicture}
\]\vskip-0.6cm
\caption{A diamond diagram, its reading together with the associated loop
  configuration and arc diagram}
\label{diagrams}
\end{figure}

A \defn{diamond diagram} of rank-$N$ is a trapezoidal arrangement of boxes
with $N-1$ rows starting with a north-east diagonal and ending with south-east
diagonal, where each box can be either black or white. The rows are indexed
from bottom to top. The \defn{reading} of such a diagram is the sequence
$\underline{\mathrm{\bf i}}= (i_1, \dots, i_\ell)$ obtained by recording the
row index $i_k$ of the $k$-th black box, starting on the left and reading each
south-east diagonal from top to bottom. Associated to the reading
$\underline{\mathrm{\bf i}}$ is the monomial
$\E{\underline{\mathrm{\bf i}}} := \E{i_1} \cdots \E{i_\ell}$ in the Okada
algebra $\Okada(X,Y)$.  We identify diagrams differing by empty south-east
diagonals on the right.  This identification is compatible with the reading
and the associated monomial. See~\cref{diagrams} for some examples. Using
rewriting techniques on such diagrams, one can show right away that the Okada
algebra has dimension $N!$.

The relevant combinatorics becomes more transparent
after we re-encode a diamond diagram as a
 \defn{fully packed loop configuration} (FPLC). This is done by
replacing black and white squares respectively with double U-turn and double
horizontal squares:
\mbox{$\begin{tikzpicture}[diamond, diamond/bound={1}{1}
    ]\matrix[diamond/matrix] {
      \bl\\
    };
\end{tikzpicture}\to
\begin{tikzpicture}[diamond, diamond/bound={1}{1}
  ]\matrix[diamond/matrix] {
    \tu\\
};
\end{tikzpicture}$}
and
\mbox{$\begin{tikzpicture}[diamond, diamond/bound={1}{1}
    ]\matrix[diamond/matrix] {
    \ \\
};
\end{tikzpicture}\to
\begin{tikzpicture}[diamond, diamond/bound={1}{1}
  ]\matrix[diamond/matrix] {
    \st\\
};
\end{tikzpicture}$}. The paths fragments at the top and bottom of the
trapezoid are completed by adding horizontal lines. The result is a set
$\mathcal{C}$ of non-crossing planar loops and arcs. The endpoints of the arcs
are situated on the left and right boundaries of the trapezoid and we number
these endpoints, from bottom to top, with positive indices (on the left) and
negative indices (overlined, on the right).  See \cref{diagrams} where we've
colored some of the arcs in order to make the picture more legible.

\renewcommand{\intextsep}{0cm}
The horizontal arc segments  in the
\mbox{$\begin{tikzpicture}[diamond, diamond/bound={1}{1}
]\matrix[diamond/matrix] {
    \st\\
};
\end{tikzpicture} \, \text{-boxes}$}
occupy levels ${1, \dots, N}$ starting from the bottom
of the trapezoid. The \defn{height} of an arc/loop in an FPLC is the minimal level
of the horizontal segments which form it.
The height statistic of an arc/loop is 
invariant under the following local moves which 
implement the Okada relations:
\vskip-0.4cm

\[
\OSys\eqdef
\left\{
\begin{tikzpicture}[diamond,
  diamond/bound={2}{2}]\matrix[diamond/matrix] {
    \nn \\
    & \st & \st \\
    & \st & \tu & \nn& \nn\\
  };
\end{tikzpicture}
\longmapsto
\begin{tikzpicture}[diamond,
  diamond/bound={2}{2}]\matrix[diamond/matrix] {
    \nn \\
    & \tu & \st \\
    & \st & \st & \nn& \nn\\
  };
\end{tikzpicture}
\ ,\quad
\begin{tikzpicture}[diamond,
  diamond/bound={2}{2}]\matrix[diamond/matrix] {
    \nn \\
    & \tu & \st \\
    & \st & \tu & \nn& \nn\\
  };
\end{tikzpicture}
\longmapsto\ x_i\,
\begin{tikzpicture}[diamond,
  diamond/bound={2}{2}]\matrix[diamond/matrix] {
    \nn \\
    & \tu & \st \\
    & \st & \st & \nn& \nn\\
  };
\end{tikzpicture}
\ ,\quad
\begin{tikzpicture}[diamond,
  diamond/bound={2}{2}]\matrix[diamond/matrix] {
    \nn \\
    & \tu & \tu \\
    & \st & \tu & \nn& \nn\\
  };
\end{tikzpicture}
\longmapsto\ y_{i-1}\,
\begin{tikzpicture}[diamond,
  diamond/bound={2}{2}]\matrix[diamond/matrix] {
    \nn \\
    & \tu & \st \\
    & \st & \st & \nn& \nn\\
  };
\end{tikzpicture}\right\}.
\]
The first and third moves can be viewed as restricted isotopies which
transform arcs and loops horizontally and downward, while the second move
erases loops. By repeatedly applying local moves, each FPLC can be brought to
a \defn{normal} form (\textit{ie.} a configuration without any possible
move). This normal form is independent of the sequence of moves used
to obtain it and is therefore uniquely defined. It contains no arcs which go up and then down when followed
in any direction; in particular, there are no loops. There is a bijection
between permutations and normal forms which shows the following result:
\begin{theorem}
  For any $N$ and for any specialization of the $X, Y$ parameters, the map
  $\sigma\mapsto \E{\sigma}$ is a bijection from $\SG[N]$ to the monoid
  $\Okada$ and the family $(\E{\sigma})_{\sigma\in\SG}$ is a basis for the
  Okada algebra $\Okada(X,Y)$. In particular the dimension of the Okada algebra
  $\Okada(X,Y)$ is always $N!$.
\end{theorem}

We abbreviate the structure of a FPLC $\mathcal{C}$ by removing its loops, labeling each
arc by its respective height, and taking the isotopy class of what remains. We
denote the result $[\mathcal{C}]$; an example is depicted in the third image
of \cref{diagrams}.  In view of the previous remarks
$[\mathcal{C}] = [\mathcal{D}]$ whenever $\mathcal{C}$ and $\mathcal{D}$ are
two FPLCs of rank $N$ which are related by a sequences of moves. It turns out
that this is actually an equivalence, providing us with a diagram model for the
Okada algebra which we now examine.
\bigskip

Recall that a rank $N$ \defn{non-crossing arc-diagram} is a visualization of a
perfect matching linking vertices $\{1, \dots , N\}$ and
$\{\ov{1}, \dots, \ov{N} \}$ on the right and left boundaries of a rectangle
by non-crossing arcs (drawn in the interior of the rectangle). A pair
$\{a,b\}$ in the matching is depicted by an arc joining vertices
$a,b \in [N] \cup [\ov{N}]$ and is denoted by $\arc{a}{}{b}$. Either $a,b$ are
both positive, both negative, or else have different signs; in the later case
we say the arc $\arc{a}{}{b}$ is a \defn{propagating}. Only the incidence
relations of the arc-diagram are relevant, and so isotopic diagrams are
considered equivalent.

An arc $\arc{a}{}{b}$ is said to be \defn{nested} in another arc
$\arc{c}{}{d}$ if $c<a<b<d$. Nesting defines a partial order on the arcs of a
non-crossing arc diagram. The reader should be aware that any arc situated
above a propagating arc is nested in the later. In particular, given
any pair of propagating arcs, one arc must be nested in the other;
consequently the nesting order is total when restricted to propagating
arcs.
\begin{definition}
  \label{def.okada.arc}
  A rank $N$ \defn{Okada arc-diagram} is a rank $N$ non-crossing arc-diagram
  where each arc $\arc{a}{}{b}$ is 
  assigned an \defn{height-label} ${\bf h}(\arc{a}{}{b})$ such that 
  \begin{enumerate}
  \item\label{okada-arc-cond-interv} ${\bf h}(\arc{a}{}{b})$ must
    be at least $1$ and at most $\min(|a|, |b|)$,
  \item\label{okada-arc-cond-parity} ${\bf h}(\arc{a}{}{b})$ 
    must have the same parity as $\min(|a|, |b|)$,
  \item\label{okada-arc-cond-nested} ${\bf h}(\arc{a}{}{b}) > \, {\bf h}(\arc{c}{}{d})$ 
  whenever $\arc{a}{}{b}$ is nested in $\arc{c}{}{d}$.
  \end{enumerate}
\end{definition}

The set of all Okada arc-diagrams of rank $N$ is denoted $\mathfrak{A}_N$
and $\C\mathfrak{A}_N$ will denote the vector space spanned by
all Okada arc-diagrams of rank $N$.

\begin{definition}
Given $C, D \in \mathfrak{A}_N$ their \defn{composition} $C \circ D$ 
is the diagram
obtained by merging the right boundary nodes of 
$C$ with the left boundary nodes of $D$ and
concatenating their respective arcs. 
The diagram $C \circ D$ may include loops, created from 
concatenated arc fragments of diagrams $C$ and $D$. The \defn{height label} of an arc/loop
in $C \circ D$ is the minimum of the height labels of the arc fragments of $C$ and $D$
which comprise it. Let $[C \circ D]$ denoted the isotopy class of the
height labeled, non-crossing arc-diagram obtained by removing all loops from the composition.
\end{definition}

\begin{lemma}
  If $C, D \in \mathfrak{A}_N$ then $[C \circ D] \in \mathfrak{A}_N$. Hence
  $\mathfrak{A}_N$ acquires the structure of a monoid denoted $\Okada$ whose
  unit is the \defn{identity Okada arc-diagram} $\id_N$ which consists
  entirely of labeled propagating arcs ${\bf h}(\arc{a}{}{\ov{a}}) = a$ for
  all $1\leq a \leq N$.
\end{lemma}

For simplicity we'll present the following results  
for arbitrary $X$-parameters together with 
the assumption that $y_k = 1$ for all $k \geq 1$.
This is sufficiently generic to include the semi-simple case
and all but the most degenerate specializations.
\begin{definition}
Given $C, D \in \mathfrak{A}_N$ their \defn{product} $C \cdot D$ 
is the element $X^{\underline{\ell}} \, [C \circ D]$ in $\C\mathfrak{A}_N$ where
$X^{\underline{\ell}} = \prod_{k \geq 0} x_k^{\ell_k}$ and $\ell_k$
counts the number of loops $\gamma$ in $C \circ D$ 
whose height label equals $k$.
\end{definition}

\begin{lemma}
The product $C \cdot D$ endows $\C\mathfrak{A}_N$ with the structure of 
an associative, unital algebra which we denote 
$\OkArc(X,{\bf 1})$. 
Using a rewriting rule, the $Y$-parameters can also be
incorporated in the diagram product 
and $\OkArc(X,Y)$ will denote the corresponding diagram algebra.

\end{lemma}
\cref{fig-okada-arc} below shows some examples.

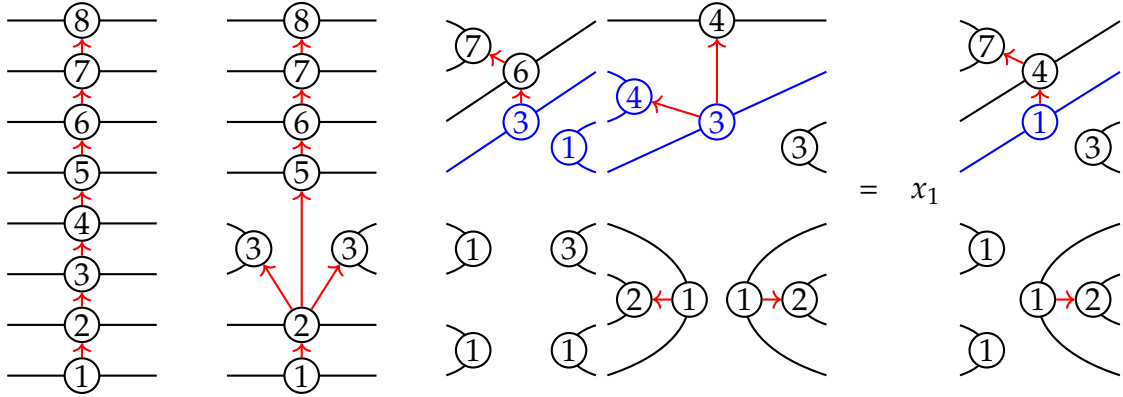
\begin{figure}[ht]
\[
  \makebox[\textwidth]{
\begin{tikzpicture}[xscale=0.7,yscale=0.45, thick, baseline={(current bounding box.east)}]
\tikzstyle{vertex} = [shape=rectangle, minimum height=15pt, minimum width=0pt, inner sep=0pt]
\tikzstyle{mid} = [draw, fill=white, shape=circle, minimum size=2pt, inner sep=1pt]
\foreach \i in {1,...,8} {
  \node[vertex] (G-\i) at ($(-1.44, \i*1.5-1.5)$) {};
  \node[vertex] (G--\i) at ($(1.44, \i*1.5-1.5)$) {};
}
\foreach \i in {1,...,8} {
   \draw (G-\i) -- (G--\i) node[pos=0.5,mid] (M-\i-\i) { $\i$ };
}
\foreach \j [count=\i from 1] in {2,...,8} {
   \draw[color=red, ->] (M-\i-\i) -- (M-\j-\j);
}
\end{tikzpicture}
\hfil
\begin{tikzpicture}[xscale=0.7,yscale=0.45, thick, baseline={(current bounding box.east)}]
\tikzstyle{vertex} = [shape=rectangle, minimum height=15pt, minimum width=0pt, inner sep=0pt]
\tikzstyle{mid} = [draw, fill=white, shape=circle, minimum size=2pt, inner sep=1pt]
\foreach \i in {1,...,8} {
  \node[vertex] (G-\i) at ($(-1.44, \i*1.5-1.5)$) {};
  \node[vertex] (G--\i) at ($(1.44, \i*1.5-1.5)$) {};
}
\draw (G--1) -- (G-1) node[pos=0.5,mid] (M-1-1) { 1 };
\draw (G--2) -- (G-2) node[pos=0.5,mid] (M-2-2) { 2 };
\draw (G-3) .. controls +(0.70, 0.35) and +(0.70, -0.35) .. (G-4) node[pos=0.5,mid] (M3-4) { 3 };
\draw (G--4) .. controls +(-0.70, -0.35) and +(-0.70, 0.35) .. (G--3) node[pos=0.5,mid] (M-4--3) { 3 };
\draw (G--5) -- (G-5) node[pos=0.5,mid] (M-5-5) { 5 };
\draw (G--6) -- (G-6) node[pos=0.5,mid] (M-6-6) { 6 };
\draw (G--7) -- (G-7) node[pos=0.5,mid] (M-7-7) { 7 };
\draw (G--8) -- (G-8) node[pos=0.5,mid] (M-8-8) { 8 };
\draw[color=red, <-] (M3-4) -- (M-2-2);
\draw[color=red, <-] (M-4--3) -- (M-2-2);
\draw[color=red, <-] (M-8-8) -- (M-7-7);
\draw[color=red, <-] (M-7-7) -- (M-6-6);
\draw[color=red, <-] (M-6-6) -- (M-5-5);
\draw[color=red, <-] (M-5-5) -- (M-2-2);
\draw[color=red, <-] (M-2-2) -- (M-1-1);
\end{tikzpicture}
\hfil
\begin{tikzpicture}[xscale=0.7,yscale=0.45, thick, baseline={(current bounding box.east)}]
\tikzstyle{vertex} = [shape=rectangle, minimum height=15pt, minimum width=0pt, inner sep=0pt]
\tikzstyle{mid} = [draw, fill=white, shape=circle, minimum size=2pt, inner sep=1pt]
\foreach \i in {1,...,8} {
  \node[vertex] (G-\i) at ($(-1.44, \i*1.5-1.5)$) {};
  \node[vertex] (G--\i) at ($(1.44, \i*1.5-1.5)$) {};
}
\draw[color=blue] (G--7) -- (G-5) node[pos=0.5,mid] (M-7-5) { 3 };
\draw (G-3) .. controls +(0.70, 0.35) and +(0.70, -0.35) .. (G-4) node[pos=0.5,mid] (M3-4) { 1 };
\draw (G--4) .. controls +(-0.70, -0.35) and +(-0.70, 0.35) .. (G--3) node[pos=0.5,mid] (M-4--3) { 3 };
\draw (G-1) .. controls +(0.70, 0.35) and +(0.70, -0.35) .. (G-2) node[pos=0.5,mid] (M1-2) { 1 };
\draw[color=blue] (G--6) .. controls +(-0.70, -0.35) and +(-0.70, 0.35) .. (G--5) node[pos=0.5,mid] (M-6--5) { 1 };
\draw (G--8) -- (G-6) node[pos=0.5,mid] (M-8-6) { 6 };
\draw (G-7) .. controls +(0.70, 0.35) and +(0.70, -0.35) .. (G-8) node[pos=0.5,mid] (M7-8) { 7 };
\draw (G--2) .. controls +(-0.70, -0.35) and +(-0.70, 0.35) .. (G--1) node[pos=0.5,mid] (M-2--1) { 1 };
\draw[color=red, <-] (M7-8) -- (M-8-6);
\draw[color=red, <-] (M-8-6) -- (M-7-5);
\end{tikzpicture}
\begin{tikzpicture}[xscale=0.7,yscale=0.45, thick, baseline={(current bounding box.east)}]
\tikzstyle{vertex} = [shape=rectangle, minimum height=15pt, minimum width=0pt, inner sep=0pt]
\tikzstyle{mid} = [draw, fill=white, shape=circle, minimum size=2pt, inner sep=1pt]
\foreach \i in {1,...,8} {
  \node[vertex] (G-\i) at ($(-2.1, \i*1.5-1.5)$) {};
  \node[vertex] (G--\i) at ($(2.1, \i*1.5-1.5)$) {};
}
\draw[color=blue] (G--7) -- (G-5) node[pos=0.5,mid] (M-7-5) { 3 };
\draw[color=blue] (G-6) .. controls +(0.70, 0.35) and +(0.70, -0.35) .. (G-7) node[pos=0.5,mid] (M6-7) { 4 };
\draw (G--4) .. controls +(-2.10, -1.05) and +(-2.10, 1.05) .. (G--1) node[pos=0.5,mid] (M-4--1) { 1 };
\draw (G-1) .. controls +(2.10, 1.05) and +(2.10, -1.05) .. (G-4) node[pos=0.5,mid] (M1-4) { 1 };
\draw (G--6) .. controls +(-0.70, -0.35) and +(-0.70, 0.35) .. (G--5) node[pos=0.5,mid] (M-6--5) { 3 };
\draw (G--8) -- (G-8) node[pos=0.5,mid] (M-8-8) { 4 };
\draw (G-2) .. controls +(0.70, 0.35) and +(0.70, -0.35) .. (G-3) node[pos=0.5,mid] (M2-3) { 2 };
\draw (G--3) .. controls +(-0.70, -0.35) and +(-0.70, 0.35) .. (G--2) node[pos=0.5,mid] (M-3--2) { 2 };
\draw[color=red, <-] (M6-7) -- (M-7-5);
\draw[color=red, <-] (M2-3) -- (M1-4);
\draw[color=red, <-] (M-3--2) -- (M-4--1);
\draw[color=red, <-] (M-8-8) -- (M-7-5);
\end{tikzpicture}\quad =\quad $x_1$\ \ 
\begin{tikzpicture}[xscale=0.7,yscale=0.45, thick, baseline={(current bounding box.east)}]
\tikzstyle{vertex} = [shape=rectangle, minimum height=15pt, minimum width=0pt, inner sep=0pt]
\tikzstyle{mid} = [draw, fill=white, shape=circle, minimum size=2pt, inner sep=1pt]
\foreach \i in {1,...,8} {
  \node[vertex] (G-\i) at ($(-1.64, \i*1.5-1.5)$) {};
  \node[vertex] (G--\i) at ($(1.44, \i*1.5-1.5)$) {};
}
\draw[color=blue] (G--7) -- (G-5) node[pos=0.5,mid] (M-7-5) { 1 };
\draw (G-3) .. controls +(0.70, 0.35) and +(0.70, -0.35) .. (G-4) node[pos=0.5,mid] (M3-4) { 1 };
\draw (G--4) .. controls +(-2.10, -1.05) and +(-2.10, 1.05) .. (G--1) node[pos=0.5,mid] (M-4--1) { 1 };
\draw (G-1) .. controls +(0.70, 0.35) and +(0.70, -0.35) .. (G-2) node[pos=0.5,mid] (M1-2) { 1 };
\draw (G--6) .. controls +(-0.70, -0.35) and +(-0.70, 0.35) .. (G--5) node[pos=0.5,mid] (M-6--5) { 3 };
\draw (G--8) -- (G-6) node[pos=0.5,mid] (M-8-6) { 4 };
\draw (G-7) .. controls +(0.70, 0.35) and +(0.70, -0.35) .. (G-8) node[pos=0.5,mid] (M7-8) { 7 };
\draw (G--3) .. controls +(-0.70, -0.35) and +(-0.70, 0.35) .. (G--2) node[pos=0.5,mid] (M-3--2) { 2 };
\draw[color=red, <-] (M7-8) -- (M-8-6);
\draw[color=red, <-] (M-3--2) -- (M-4--1);
\draw[color=red, <-] (M-8-6) -- (M-7-5);
\end{tikzpicture}}
\]\vspace{-0.5cm}
\caption{The identity, the generator $\G{3}$ and a composition of Okada
  arc-diagrams of rank $8$. The red arrows indicate the Hasse diagram of the
  nesting order.}
\label{fig-okada-arc}
\end{figure}

Let's point out a few simple remarks. The \defn{mirror} $D^\star$
of an Okada arc-diagram, obtained by reflecting $D$ horizontally,
is an Okada-arc diagram and 
the map $D \mapsto D^\star$ extends to an anti-isomorphism of $\OkArc(X,Y)$. Let $\iota_N$
denote the map from $\OkArc[N](X,Y)$ to $\OkArc[N+1](X,Y)$ adding 
the labeled, propagating arc ${\bf h}(\arc{N+1}{}{\ov{N+1}})=N+1$
to each arc-diagram. Clearly $\iota_N$ is an injective
algebra homomorphism whose image is the set of diagrams containing
${\bf h}(\arc{N+1}{}{\ov{N+1}})=N+1$.
\begin{definition}
  For $1\leq i < N$, let $\G{i}$ denote the \defn{elementary} Okada arc-diagram containing the 
  labeled arcs \ \ ${\bf h}(\arc{j}{}{\ov{j}})=j$ \, for \, $j\neq i, i+1$,\quad
  ${\bf h}(\arc{i}{}{i+1})=i$\qandq ${\bf h}(\arc{\ov{i}}{}{\ov{i+1}})=i$.
\end{definition}
The elementary Okada diagrams $\G{1}, \dots, \G{N-1}$ satisfy Okada relations
\ref{relIxy}, \ref{relCxy}, and \ref{relSxy}.
To construct the isomorphism from $\Okada(X,Y)$ to $\OkArc(X,Y)$, we first need to show that 
the elements $(\G{i})$ generate $\OkArc(X,Y)$. It is clear from
the definition that if a product ends with $\G{i}$, then the diagram contains the arc
${\bf h}(\arc{\ov{i}}{}{\ov{i+1}})=i$. The converse is actually true: If an element
$\mathrm{\bf e} \in \mathfrak{A}_N$ contains the arc ${\bf h}(\arc{\ov{i}}{}{\ov{i+1}})=i$ then it can
be factored as $\mathrm{\bf e} = \mathrm{\bf f} \cdot \G{i}$.
\begin{proposition}
  Suppose $D \in \mathfrak{A}_N$ doesn't contain
  the arc ${\bf h}(\arc{N}{}{\ov{N}})=N$. Then there exist an integer~$i$ such that $D$
  contains the arc ${\bf h}(\arc{\ov{i}}{}{\ov{i+1}})=i$. If $I$ is the largest such
  integer, then there exists a unique arc diagram $D^\flat \in \mathfrak{A}_{N-1}$ such that
   $D$ factorize as
  \[
    D=\iota_{N-1}(D^\flat) \G{N-1}\G{N-2} \cdots \G{I}\,.
  \]
\end{proposition}
By induction, this proves the following theorem:
\begin{theorem}
  The dimension of $\OkArc(X,Y)$ is $N!$ and it is generated by the elementary diagrams
  $\G{1}, \dots, \G{N-1}$. Furthermore the map sending $\E{i}$ to $\G{i}$ extends
  multiplicatively to a unique algebra isomorphism $\Theta:\Okada(X,Y) \to \OkArc(X,Y)$.
\end{theorem}
We conclude this section by making explicit the relation between fully packed loop
configurations and Okada arc-diagrams:
\begin{proposition}
  For simplicity assume $Y = \mathrm{\bf 1}$.  
  Let $\mathcal{C}$ be a FPLC of rank $N$,
  let ${\bf \underline{i}}$ be its reading,  
  and let $\E{\bf \underline{i}}$ be the corresponding monomial in the
  Okada algebra $\Okada(X,\mathrm{\bf 1})$. Then $\Theta(\E{\bf \underline{i}}) = X^{\underline{\ell}} \, [\mathcal{C}]$ 
  where $X^{\underline{\ell}} = \prod_{k \geq 1} x_k^{\ell_k}$ and $\ell_k$
  counts the number of loops in $\mathcal{C}$ with height $k$.
\end{proposition}

\section{Fomin correspondence and Okada arc-diagrams}

We have a bijection between $\SG$ and the monoid $\OkArc$ of Okada
arc-diagrams, however, it is circuitous: Starting from a permutation, first
its code is computed, then the associated FPLC is drawn, from which an Okada
arc-diagram is obtained. It is not obvious, for example, that the inverse of a
permutation corresponds to the mirror of the associated Okada arc-diagram. The
goal of this section is to better explicate this graphical bijection which
turns out to be an incarnation of Fomin's Robinson-Schested correspondence for
the Young-Fibonacci~\cite{Stanley1988,Fomin1995} lattice. Recall that this is
a bijection between permutations of $\SG[N]$ and pairs of saturated chains in
the Young-Fibonacci lattice starting at $\emptyset$ and sharing a common
endpoint in $\YF_N$.  The reader who is not familiar with Fomin's construction
should refer to~\cite{Fomin1995}. See~\cref{fig-arc-chain} for an
example.  We will see in this section that Okada arc-diagrams are also in
natural bijection with the same pairs of chains. 

Cutting a labeled arc-diagram $D$ in the middle gives a
natural notion of a \defn{Okada half arc-diagram} containing either a labeled full arc
${\bf h}(\arc{a}{}{b})$ joining two nodes $a, b \in\setN$, or else a labeled half arc
${\bf h}(\harc{a}{})$ with a free end. Such a half arc is called \defn{propagating}.
The \defn{bra} $\langle D |$ is the Okada half arc-diagram obtained by restricting $D$
to its positive part. The \defn{ket} $| D \rangle$ is defined to be
the bra $\langle D^\star |$ of the mirror $D^\star$.

\begin{definition}
  The \defn{propagating label set} of an Okada half arc-diagram $H$ is the subset
  $\PropLab(H)$ of $\setN$ consisting of the height labels of its propagating arcs.
\end{definition}
The propagating label set of an Okada half arc-diagram of rank $N$ is always a
Fibonacci set of rank $N$.
The following trivial lemma-definition is of great importance:
\begin{lemma}[Gluing lemma]\label{lemma:gluing}
  For any Okada arc diagram $D$, the left and right half diagram $\leftd{D}$
  and $\rightd{D}$ are two Okada half arc diagrams which have the same
  propagating labels set. As consequence, it makes sense to define
  $\PropLab(D)\eqdef\PropLab(\leftd{D})=\PropLab(\rightd{D})$.

  Moreover if $L$ and $R$ are two Okada half arc-diagrams such that
  $\PropLab(L)=\PropLab(R)$, there is a unique Okada
  arc-diagram~\defn{$\glue{L}{R}$} such that $\leftd{\glue{L}{R}}=L$ and
  $\rightd{\glue{L}{R}}=R$.
\end{lemma}
To convert half arc diagrams to chains we need to restrict the former:
\begin{definition}
  For $r\leq N$, the \defn{$r$-restriction} of an Okada half arc diagram $H$ is
  the Okada half arc-diagram denoted by $H/\setN[r]$ of rank $r$ possessing
  \begin{itemize}
  \item a full arc ${\bf h}(\arc{a}{}{b})=h$ whenever $H$ contains a full arc
    ${\bf h}(\arc{a}{}{b})=h$ with $a, b \leq r$
  \item a half arc ${\bf h} (\harc{a}{})=h$ whenever $H$ contains either a full arc 
    ${\bf h}(\arc{a}{}{b}) = h$ with $a \leq r<b$ or a half arc ${\bf h}(\harc{a}{})=h$ with
    $a \leq r$.
  \end{itemize}
\end{definition}
If $r\leq s$ then clearly the $r$-restriction of the $s$-restriction of any
Okada half arc-diagram $H$ coincides with the $r$-restriction of $H$.
\begin{definition}
  To any Okada half arc-diagram $H$ of rank $N$ we associate the sequence of
  Fibonacci sets $\defn{$\Chain(D)$}\eqdef(C_0,\dots,C_N)$ defined by
  $C_i\eqdef\PropLab(H/\setN[i])$.
\end{definition}
\begin{proposition}\label{prop:chain-diag-bij}
  The map $\Chain$ is a bijection between Okada half arc-diagrams and
  saturated chains of rank $N$ in the $\YFS$-lattice. See
  \cref{fig-arc-chain} for an example.
\end{proposition}
\begin{theorem}
  Given a permutation $\sigma \in \SG[N]$, let $L_\sigma$ and $R_\sigma$
  denote the two Okada half arc-diagrams associated to the pair of saturated
  chains obtained from Fomin's RS-correspondence.  Then
  $\Theta(\E{\sigma})= \glue{L_\sigma}{R_\sigma}$. Moreover
  $\Theta(\E{\sigma})^\star = \Theta(\E{\sigma^{-1}}) =
  \glue{R_\sigma}{L_\sigma}$.
\end{theorem}
\section{Structure of the Okada algebra and monoid}

From now on, we identify $\Okada(X,Y)$ and $\OkArc(X,Y)$ through the isomorphism
$\Theta$. The goal of this section is to understand the structure
of the Okada algebra and its monoid via the $\YFS_N$ \defn{dominance
  order}. In particular, we describe the stratification of the Okada algebra
by two-sided ideals (generated by \defn{free elements}) and the Green
relations for the monoid (which allows us show that the monoid is aperiodic).
This allows us to prove cellularity of the Okada algebra in the next section.
\medskip

\begin{definition}
  Let $S=\{s_1<\dots<s_k\}$ and $T=\{t_1<\dots<t_\ell\}$ be two Fibonacci sets of
  the same rank $N$. We says that \defn{$S$ is dominated by $T$} and write
  $S \isdomeq T$ if $k<\ell$ and $s_{k-i}\leq t_{\ell-i}$ for any $0\leq i<k$. We
  write $S \isdom T$ if $S \isdomeq T$ but $S \neq T$.
\end{definition}
\begin{proposition}
  $(\YFS_N, \isdomeq)$ is a ranked lattice.
\end{proposition}

\begin{figure}[ht]
  \begin{subfigure}{.475\linewidth}
    \def\es{\emptyset}
      \begin{tikzpicture}[xscale=0.70,yscale=0.5, thick, baseline={(current bounding box.east)}]
\tikzstyle{vertex} = [shape=circle, minimum size=7pt, inner sep=1pt]
\tikzstyle{mid} = [draw, fill=white, shape=circle, minimum size=2pt, inner sep=1pt]
\foreach \i in {0,...,8} {
  \node[vertex] (G-\i) at ($(-2.64, \i*1.5-1.5)$) {};
  \node[vertex] (G--\i) at ($(2.64, \i*1.5-1.5)$) {};
}
\draw (G--3) -- (G-1) node[pos=0.5,mid] (M-3-1) { 1 };
\draw (G--8) .. controls +(-2.592, -1.296) and +(2.592, 1.296) .. (G-4) node[pos=0.5,mid] (M-8-4) { 2 };
\draw (G-5) .. controls +(2.1, 1.05) and +(2.1, -1.05) .. (G-8) node[pos=0.5,mid] (M5-8) { 3 };
\draw (G--2) .. controls +(-0.7, -0.35) and +(-0.7, 0.35) .. (G--1) node[pos=0.5,mid] (M-2--1) { 1 };
\draw (G-2) .. controls +(0.7, 0.35) and +(0.7, -0.35) .. (G-3) node[pos=0.5,mid] (M2-3) { 2 };
\draw (G--7) .. controls +(-2.1, -1.05) and +(-2.1, 1.05) .. (G--4) node[pos=0.5,mid] (M-7--4) { 4 };
\draw (G-6) .. controls +(0.7, 0.35) and +(0.7, -0.35) .. (G-7) node[pos=0.5,mid] (M6-7) { 6 };
\draw (G--6) .. controls +(-0.7, -0.35) and +(-0.7, 0.35) .. (G--5) node[pos=0.5,mid] (M-6--5) { 5 };
\draw[color=red, <-] (M6-7) -- (M5-8);
\draw[color=red, <-] (M5-8) -- (M-8-4);
\draw[color=red, <-] (M2-3) -- (M-3-1);
\draw[color=red, <-] (M-6--5) -- (M-7--4);
\draw[color=red, <-] (M-7--4) -- (M-3-1);
\draw[color=red, <-] (M-8-4) -- (M-3-1);
\foreach \x [count=\xi from 0] in {
  \es, \{1\}, {\{1,2\}}, \{1\}, {\{1,2\}}, {\{1,2,3\}}, {\{1,2,3,6\}},
  {\{1,2,3\}}, {\{1,2\}}
} {
  \node [left] at (G-\xi) {$\x$};
}
\foreach \x [count=\xi from 0] in {
  \es, \{1\}, \es, \{1\}, {\{1,4\}}, {\{1,4,5\}}, {\{1,4\}}, {\{1\}},
  {\{1,2\}}
} {
  \node [right] at (G--\xi) {$\x$};
} 
\end{tikzpicture}
\caption{An arc-diagram with its associated chains}
\label{fig-arc-chain}
\end{subfigure}
%
\begin{subfigure}{.475\linewidth}
  \[\begin{tikzpicture}[>=latex,line join=bevel,draw,black,
      baseline={(current bounding box.east)}]
\node (50) at (6.5,0) {$\left\{1\right\}$};
\node (51) at (6.5,1) {$\left\{3\right\}$};
\node (52) at (5.5,2) {$\left\{1, 2, 3\right\}$};
\node (53) at (7.5,2) {$\left\{5\right\}$};
\node (54) at (6.5,3) {$\left\{1, 2, 5\right\}$};
\node (55) at (6.5,4) {$\left\{1, 4, 5\right\}$};
\node (56) at (6.5,5) {$\left\{3, 4, 5\right\}$};
\node (57) at (6.5,6) {$\left\{1, 2, 3, 4, 5\right\}$};
\draw [black,->] (50) -- (51);
\draw [black,->] (51) -- (52);
\draw [black,->] (51) -- (53);
\draw [black,->] (52) -- (54);
\draw [black,->] (53) -- (54);
\draw [black,->] (54) -- (55);
\draw [black,->] (55) -- (56);
\draw [black,->] (56) -- (57);
\end{tikzpicture}\]
\caption[Dominance orders of Fibonacci sets]{The dominance order on $\YFS_5$}
\label{fig.dominance}
\end{subfigure}
\end{figure}

\begin{definition}
A \defn{free set} of rank $N$ is a subset of $[N]$
which does not contain both $i$ and $i+1$ for all $1 \leq i < N$. 
The map
$S \mapsto \mathfrak{F}(S) := \{ i \mid i - \max\{ s \in S \mid s < i \} \ \text{is odd} \}$
defines a bijection from the collection of rank $N$ Fibonacci sets to the collection
of rank $N$ free sets. 
\end{definition}

\begin{definition}
For $S \in \YFS_N$ the corresponding \defn{free element} in
$\Okada(X,Y)$ is $\E{S} := \prod_{i \in \mathfrak{F}(S)} \E{i}$. 
\end{definition}
Note that $\E{S} = \E{\sigma_S}$ where $\sigma_S := \prod_{i \in \mathfrak{F}(S)} \sigma_i$
is the associated \defn{free involution} in $\SG[N]$. 

\begin{remark}
The half arc-diagrams 
$\langle \E{S} |$ and $| \E{S} \rangle$
coincide for any $S \in \YFS_N$.
Furthermore $\langle \E{S} |$
consists only of
labeled propagating arcs ${\bf h}(\harc{s}{})= s$ 
for $s \in S$ and
labeled full arcs ${\bf h}(\arc{i}{}{i+1})= i$
for $i \in \mathfrak{F}(S)$. 
\end{remark}

\begin{proposition}
Let  $\mathfrak{J}_S$ be the two-sided ideal in 
$\Okada(X,Y)$ generated a free element $\E{S}$
for $S \in \YFS_N$, then
$\mathfrak{J}_S \subseteq \mathfrak{J}_T$ if and only if $T \preceq S$
for any pair $S, T \in \YFS_N$.
\end{proposition}

\begin{theorem}[Triangular Factorization]\label{theorem:triangular-factor}
  For $\sigma \in \SG[N]$ there exists a unique pair of permutations $\rho, \tau \in \SG[N]$ such that 
  $\E{\sigma} = \E{\rho} \cdot \E{S} \cdot \E{\tau}$ where
  $\ell(\sigma) = \card S + \ell(\rho) + \ell(\tau)$
  and $S \! \preceq \inf ( \PropLab (\E{\rho} ), \PropLab (\E{\tau} ) )$
  and where $S = \PropLab (\E{\sigma})$. 
\end{theorem}

Returning to the Okada monoid,
an element $\mathrm{\bf e}\in\Okada$ is said to be \defn{involutive} whenever it equals its mirror $\mathrm{\bf e}^\star$.
Involutive elements are always idempotents and 
thanks to the $\RS$ correspondence, these are precisely the basis monomials $\E{\sigma}$
where $\sigma \in \SG[N]$ is an involution (\ie $\sigma^2=1$).
\begin{remark}
  The set of idempotents in $\Okada$ is not exhausted by the involutive elements.
  For example in $\Okada[3]$ all element are idempotents, while
  $\E{1}\E{2}\E{3}$ and $\E{3}\E{2}\E{1}$ are the only non-idempotents in
  $\Okada[4]$. Computer calculations show that the number of idempotents for
  $N \leq 10$ are: $1,\ 1,\  2,\  6,\  22,\  108,\  594,\  4116,\  30500,\  274006,\  2560400.$
\end{remark}

\begin{proposition}\label{lemma:left-eq-or-prop-le}
  Let $\mathrm{\bf e,f}\in\Okada$. Then either $\leftd{\mathrm{\bf ef}}=\leftd{\mathrm{\bf e}}$ and thus
  $\PropLab(\mathrm{\bf ef})=\PropLab(\mathrm{\bf e})$ or $\PropLab(\mathrm{\bf ef})\isdom\PropLab(\mathrm{\bf e})$. As a
  consequence, $\PropLab(\mathrm{\bf ef})\isdomeq\inf(\PropLab(\mathrm{\bf e}), \PropLab(\mathrm{\bf f}))$.
\end{proposition}
The previous proposition is the main ingredient of the following theorem which
describe the structure of the Okada monoid:
\begin{theorem}
  The monoid $\Okada$ is aperiodic, i.e. there exists an integer $K$
  such that $\mathrm{\bf e}^K=\mathrm{\bf e}^{K+1}$ for all $\mathrm{\bf
    e}\in\Okada$. Equivalently, all the groups in $\Okada$ are trivial.
\end{theorem}
Recall that $\RR$ (resp. $\JJ$) is the equivalence relation on $\Okada$ such
that $\mathrm{\bf e} \! \RR \! \mathrm{\bf f}$ if $\mathrm{\bf e}$ and
$\mathrm{\bf f}$ generate the same right (resp. two-sided) ideals.
\begin{theorem}\label{prop:r-j-class}
  Each $\RR$-class of $\Okada$ contains a unique involutive element. 
  Each $\JJ$-class of $\Okada$ contains a unique free element. Moreover, the
  free representative of $\mathrm{\bf e}\in\Okada$ is the free element having the same
  propagating set as $\mathrm{\bf e}$.
\end{theorem}
\section{Cellular structure of the Okada algebra}

Recall that a cellular algebra $A$ is a finite dimensional algebra with
distinguished cellular basis which is particularly well-adapted to studying
the representation theory of $A$, especially as the ground ring/field 
varies. For brevity, we skip a general discussion
about cellular algebras and point the reader to~\cite{Graham1996CellularA} for
definitions and context.

\begin{definition}
Let $\mathfrak{H}_N$ and $\C\mathfrak{H}_N$ denote respectively the set 
and the vector space spanned by all Okada half arc-diagrams of rank $N$.
Likewise $\mathfrak{H}_N^S$ and $\C\mathfrak{H}_N^S$
will denote the set and the vector space spanned by all half diagrams $H \in \mathfrak{H}_N$
for which $\PropLab(H)=S$ where $S \in \YFS_N$.
We extend the bra map $D \mapsto \langle D |$
by linearity to obtain a map from $\Okada(X,Y)$ to $\C\mathfrak{H}_N$.
\end{definition}

The following result is a consequence of the factorization given in Proposition \ref{theorem:triangular-factor}:
\begin{theorem}
  The Okada algebra $\Okada(X, Y)$ is cellular with the following data
  \begin{enumerate}
  \item\label{cond-cell-poset} A cell-poset is $\Lambda_N =(\YFS_N, \isdomeq)$.
  \item\label{cond-cell-invar} An index set $M_S = \mathfrak{H}_N^S$ for each $S \in \YFS_N$
  \item \label{cond-cell-basis} A cellular basis element $C^S_{L, R}\eqdef \glue{L}{R}$
  associated to $L, R\in \mathfrak{H}_N^S$
  \item\label{cond-cell-invol} An involutive anti-isomorphism given by
    the mirror map $\star: \glue{L}{R} \mapsto \glue{R}{L}$. 
  \end{enumerate}
\end{theorem}

\begin{remark}
  The left $\Okada(X,Y)$ \defn{cell module} associated to $S \in \YFS_N$ can
  be realized by the vector space $\C\mathfrak{H}_N^S$ equipped with the left
  action defined by
\[ D  \bullet H \ := \begin{cases}
\big\langle D \cdot (\glue{H}{H}) \big|
&\text{if $\PropLab\big(D \circ (\glue{H}{H})\big) = S$} \\
0
&\text{otherwise}
\end{cases}
\]
where $D \in \Okada(X,Y)$ and $H \in \mathfrak{H}_N^S$. For generic values of the $X,Y$ parameters $\C\mathfrak{H}_N^S$
is irreducible.\footnote{The cell module $\mathfrak{H}_N^S$ carries an invariant bilinear form $\varphi_S$. 
We conjecture an explicit value for the determinant of the associated Gram matrix $G_S$, which we
express in terms of (specialized) clone Schur functions.}
\end{remark}

\section{Prospectives}
For a fixed choice of a threshold $k \geq 1$, we \defn{truncate} any Okada
arc-diagram $D$, replacing its height labels $h$ by $\min(h, k)$. The
$k$-truncated Okada arc-diagrams form a multiplicative basis for a
\defn{higher Blob algebra} $\mathrm{Blob}^{(k)}_N$, which can be realized as a
quotient of the Okada algebra $\Okada(X,Y)$ after specializing the $X,Y$
parameters appropriately. In particular the Temperely-Lieb and Martin-Saleur
Blob algebras~\cite{Martin1994a} are recovered for $k = 1,2$ respectively.  It
seems that the corresponding Bratelli diagram $\YF^{(k)}$ naturally embeds
into the $\YF$-lattice and can be seen as a Fibonacci counterpart of the
sublattice of integer partitions with at most $k$ parts. Both the
Temperely-Lieb and the Blob algebras are {\it intertwiner algebras} which
raises the question of whether the higher Blob algebras have such a
description for $k \geq 3$. If so, this would be indicative of a Fibonacci
version of {\it Schur-Weyl duality}, and would entail, on a combinatorial level,
a well-behaved version of the RSK-correspondence.

It should be possible to incorporate height labels into other
diagram algebras such as the partition and Brauer algebras.
Can one, for example, define a suitable notion of height labeled {\it braids}
together with {\it skein relations} consistent with
these labels? A satisfactory answer might shed light onto
the problem of identifying appropriate {\it Jucys-Murphy} elements for the Okada
algebras.

\printbibliography[heading=bibintoc] 
\label{sec:biblio}

\end{document}